\documentclass[10pt,reqno]{amsart}
\usepackage[utf8]{inputenc} 
\usepackage[T1]{fontenc} 
\usepackage[french,english]{babel} 
\usepackage{lmodern} 
\usepackage{enumerate} 
\usepackage{geometry} 
\usepackage{fancyhdr} 
\usepackage{verbatim} 
\usepackage{amsmath,amssymb, amsthm} 
\usepackage{dsfont} 
\usepackage{graphicx} 
\usepackage{tikz} 
\usepackage{float} 
\usetikzlibrary{matrix,arrows}
\usepackage{url}
\usepackage{xypic} 
\usepackage[all]{xy}
\usepackage{stmaryrd}

\usepackage{hyperref} 
\hypersetup{pdfstartview=XYZ} 

\usepackage[final]{pdfpages} 

\usepackage[leqno]{amsmath}     
\usepackage{amssymb}            
\usepackage{mathrsfs}           
\usepackage{dsfont}             
\usepackage{amsxtra}            

\newtheoremstyle{th}{8pt}{8pt}{\itshape}{}{\bfseries}{ ---}{5pt}{\thmname{#1}\thmnumber{ #2} \thmnote{\bfseries (#3).}}
\theoremstyle{th}

\theoremstyle{break}
\newtheorem{Theorem}{Theorem}[section]

\newtheorem{Lemma}[Theorem]{Lemma}

\newenvironment{Proof}{\begin{proof}} {\end{proof}}

\newtheoremstyle{def}{8pt}{8pt}{}{}{\bfseries}{}{5pt}{\thmname{#1}\thmnumber{ #2} \thmnote{\bfseries (#3).}}
\theoremstyle{def}
\newtheorem{Definition}[Theorem]{Definition}

\newtheorem{Remark}[Theorem]{Remark}

\newtheorem{Notation}[Theorem]{Notation}

\theoremstyle{nonumberplain}

\newcommand{\field}[1]{\mathds{#1}}
\newcommand{\Z}{\field{Z}}              
\newcommand{\N}{\field{N}}              

\newcommand{\C}{\field{C}}              

\DeclareMathOperator{\Spec}{Spec}

\DeclareMathOperator{\car}{char}

\DeclareMathOperator{\CH}{CH}

\DeclareMathOperator{\Ker}{ker}

\DeclareMathOperator{\van}{van}

\DeclareMathOperator{\Alb}{Alb}

\usepackage[scaled=.90]{helvet}

\begin{document}

\date{}
 
\title{\large{The kernel of the Gysin homomorphism for \; \\  smooth projective curves}}

\author{C. Schoemann}\thanks{The author was supported by the project “Group schemes, root systems, and
related representations” founded by the European Union - NextGenerationEU through
Romania’s National Recovery and Resilience Plan (PNRR) call no. PNRR-III-C9-2023-
I8, Project CF159/31.07.2023, and coordinated by the Ministry of Research, Innovation
and Digitalization (MCID) of Romania.}

\maketitle


\vspace{-0.7cm}

\begin{center}
For my friend Dr. Rina Paucar
\end{center}

\begin{abstract}
Let $S$ be a smooth projective connected surface over an algebraically closed field $k$ and $\Sigma$ the linear system of a very ample divisor $D$ on $S$.
Let $d:=\dim(\Sigma)$ be the dimension of $\Sigma$ and 
$\phi_{\Sigma}: S \hookrightarrow \mathbb{P}^{d}$
the closed embedding of $S$ into $\mathbb{P}^{d}$,
induced by $\Sigma$.
For any closed point $t\in\Sigma \cong\mathbb{P}^{d^*}$, let $C_t$ be the corresponding hyperplane section on $S$, and let $  r_t:C_t\hookrightarrow S$
be the closed embedding of the curve $C_t$ into $S$.
Let $\Delta:= \{t \in \Sigma: C_t  \text{ is singular}\}$ be the discriminant locus of $\Sigma$ and let $U :=\Sigma\setminus \Delta$.

For $t \in U$, the kernel of the Gysin homomorphism of the Chow groups of $0$-cycles of degree zero, from
$\CH_0(C_t)_{\deg=0}$ to $\CH_0(S)_{\deg=0}$
is the countable union of shifts of a certain abelian subvariety $A_t$ inside $J(C_t)$, the
Jacobian of the curve $C_t$ (\cite{PS24} for $k \cong \C$, \cite{SW25} for $k \cong \overline{\mathbb{F}_q((t))}$). We prove that for every closed point $t \in U$
either
$A_t$ coincides with the abelian variety $B_t$ inside $J(C_t)$ corresponding to the vanishing cohomology $H^1(C_t, k')_{\van}$, where $k'$ is the minimal field of definition of $k$, and then the Gysin kernel is a countable union of shifts of $B_t$,  or $A_t = 0$, in which case the Gysin kernel
is countable. Using the language of algebraic stacks as a generalisation of algebraic varieties this is done by constructing an increasing filtration of Zariski countable open substacks $U_i, i \in I,$ of $U$, where $I$ is a countable set and by applying a convergence argument.
\end{abstract}

\section{Introduction}
Let $S$ be a smooth projective connected surface over an algebraically closed field $k$ and $\Sigma$ the linear system of a very ample divisor $D$ on $S$.
Let $d:=\dim(\Sigma)$ be the dimension of $\Sigma$ and 
$\phi_{\Sigma}: S \hookrightarrow \mathbb{P}^{d}$
the closed embedding of $S$ into $\mathbb{P}^{d}$,
induced by $\Sigma$.
For any closed point $t\in\Sigma \cong\mathbb{P}^{d^*}$, let $C_t$ be the corresponding hyperplane section on $S$, and let $  r_t:C_t\hookrightarrow S$
be the closed embedding of the curve $C_t$ into $S$.
Let $\Delta$ be the discriminant locus of $\Sigma$, that is, $    \Delta :=\{t\in \Sigma: C_t\text{ is singular}\}$.
Then $U :=\Sigma\setminus \Delta = \{t\in \Sigma: C_t\text{ is smooth}\}$. Let $k'$ be the minimal field of definition of $k$ (\cite{Wei62}) and $o_{k'}$ its ring of integers. Let
\begin{equation*}
    r_{t}^*:H^{1}(C_t,\mathbb{Z}_l)\rightarrow H^3(S,\mathbb{Z}_l)
\end{equation*}
be the Gysin homomorphism on étale cohomology groups induced by $r_t$,  whose kernel $H^{1}(C_t,\mathbb{Z}_l)_{\van}$ is called the
\emph{vanishing cohomology of} $C_t$ (\cite{Voi02}, 3.2.3).

\begin{Remark}
 If $\text{char}(k)=p$, then we chose $\mathbb{Z}_l$ with the prime $l$ such that $(l,p) = 1$.
\end{Remark}

\begin{Remark}
If $k$ is the algebraic closure of a non-archimedean local field $k'$ at a prime $p$, we choose the coefficients of the cohomology groups to be the ring of integers $o_l$ where $l$ is a prime power such that $l$ and $p$ are coprime.
\end{Remark}

Let $J_t=J(C_t)$ 
be the Jacobian of the curve $C_t$ and let $B_t$ be the abelian subvariety of the abelian variety $J_t$ corresponding to the
substructure $H^{1}(C_t,\mathbb{Z}_l)_{\van}$ of $H^{1}(C_t,\mathbb{Z}_l)$. 

\smallskip

Let $\CH_0(S)_{\deg=0}$ be the Chow group of $0$-cycles of degree zero on $S$, and for any closed point $t\in\Sigma$, let
$\CH_0(C_t)_{\deg=0}$ be the Chow group of $0$-cycles of degree zero on $C_t$.

\smallskip

For any closed point $t\in\Sigma$, let
\begin{equation*}
    r_{t}^*:\CH_0(C_t)_{\deg=0}\rightarrow \CH_0(S)_{\deg=0}
\end{equation*}
be the Gysin pushforward homomorphism on the Chow groups of degree zero $0$-cycles of $C_t$ and $S$, respectively,
induced by $r_t$,  whose kernel
\begin{equation*}
    G_t=\Ker(r_{t}^*)
\end{equation*}
is called the \emph{Gysin kernel} associated with the hyperplane section $C_t$. 

\bigskip


Let $U = \Sigma \setminus \Delta$. We have the following

\begin{Theorem}\label{Pau22}{\cite{PS22}, \cite{PS24}}
Let $k \cong \C$.
\begin{itemize}
    \item[(a)] For each $t\in U$ there is an abelian variety $A_t\subset B_t$ such that
    \begin{equation*}
        G_t=\Ker(r_t^*)=\bigcup_{\text{countable}}\text{translates of} \; A_t
    \end{equation*}
    \item[(b)] For a very general $t \in U$ (i.e. for every $t$ in a c-open subset $U_0$ of $U$) either
    \begin{itemize}
 \item[1.] $A_t = B_t$, and then $G_t = \bigcup_{\text{countable}} \text{translates of} \; B_t$, or
 \item[2.] $A_t = 0$, and then $G_t$ is countable.
    \end{itemize}
    \item[(c)] If $alb_S: \CH_0(S)_{\hom}\rightarrow \Alb(S)$ is not an isomorphism,  for a very general $t$ in $U$, 
    then $G_t$ is countable.
\end{itemize}
\end{Theorem}

The subset $U_0$ being countable open in \ref{Pau22} (b) allows to apply for all $t$ in $U_0$ in a uniform way the irreducibility of the monodromy representation on the vanishing cohomology of
a smooth section (see \cite{DK73}, \cite{D74} for the étale cohomology, \cite{La81} for the singular cohomology and \cite{Voi02}
in a Hodge theoretical context for complex algebraic varieties).
This is done by viewing $U = \Sigma \setminus \Delta$ as an integral algebraic scheme over $k$ and by passing to the general fiber,
i.e. for each closed point $t$ in $U_0$ there exists a scheme-theoretic isomorphism to the geometric generic point
$ \overline{\xi}$ over $k'$, where $k'$ is the minimal field of definition of $S$ \cite{Wei62}. This induces a scheme-theoretic isomorphism $ \kappa_t$ between the corresponding varieties $C_t$ and
$C_{\overline{\xi}}$ over $k'$ which induces an isomorphism
$ \kappa'_t $ between $A_t$ and
$A_{\overline{\xi}}$ compatible with the isomorphism on Chow groups induced by the isomorphism $ \kappa_t$. Then by \cite{BG20}
$ \kappa'_t(A_t) = A_{\overline{\xi}}$ and
$ \kappa'_t(B_t) = B_{\overline{\xi}}$ for every $k$-point in $U_0$.

\bigskip

A smooth scheme is a generalisation of the concept of a non-singular algebraic variety. Using
the language of algebraic stacks in this paper we prove the following
\begin{Theorem}{\label{PS23}}
Let $U:= \Sigma \setminus \Delta$. There exists a convergent stratification $\{\mathcal{U}_{\alpha}\}_{\alpha \in A}$ of $U$ by countable open substacks for each of which the irreducibility of the monodromy representation applies in a uniform way for all $t \in \mathcal{U}_{\alpha}$ such that the monodromy argument applies in a uniform way for all $t \in U$, seen as the set-theoretic directed union $U = \underset{\underset{i}{\rightarrow}}{\cup} \; \mathcal{U}_i$.
\end{Theorem}

In \cite{SW25} the results \ref{Pau22} $(a)$ and $(b)$ were extended to an algebraically closed uncountable field of positive characteristic. As up to (non-canonical) isomorphism there is a unique algebraically closed field of characteristic $p$ of given cardinality $\kappa$ we may think of $\overline{\mathbb{F}_q((t))}$.

Here we describe the Gysin kernel $G_t$ for the points $t$ in $U \setminus U_0$ by the local and global monodromy representations,
i.e. the action of the fundamental groups $ \pi_1(V,t)$,
where $ V = (\Sigma \setminus \Delta) \cap D$ with
$D$ a line containing $t$ in the dual space $\mathbb{P}^{d^*}$ such that $f_D$ is a Lefschetz pencil for $S$,
and $ \pi_1(U, t)$ on the vanishing cohomology $ H^1(C_t, k')_{\van}$ for an algebraically closed field $k$ of zero or of positive characteristic.
The approach is to consider algebraic stacks as a generalisation of 
algebraic varieties and to construct a stratification $\{\mathcal{U}_i \subseteq U \}_{i \in I}$ of $U$ by countable open
substacks for
each of which the monodromy argument applies for all $t \in \mathcal{U}_i$ in a uniform way (i.e. for a partially ordered, at most countable set $I$ we have
$ \mathcal{U}_i \subseteq \mathcal{U}_j$ if $ i \leq j$ and two additional conditions on $I$ for finiteness from below
and for every map $ \alpha: \Spec(k) \rightarrow U$ the set $ \{i \in I: \alpha \; \text{factors through} \; \mathcal{U}_i \}$ has a smallest
element (cf. \cite{GL19}, Def. 5.2.1.1). 

We then apply a convergence argument for the stratification $ \{\mathcal{U}_i\}_{i \in I}$ 
(cf. \cite{GL19}, Def. 5.2.2.1) 
such that the monodromy argument applies in a uniform way for all $t \in U$, seen as the set-theoretic directed union 
$U = \underset{\underset{i}{\rightarrow}}{\cup} \; \mathcal{U}_i$

\bigskip

We will begin by reviewing some background on algebraic stacks.

\section{Algebraic Stacks}

\subsection{Affine schemes, schemes and algebraic stacks}

\begin{Definition}[Ringed space]
 A ringed space $(X, \mathcal{O}_X)$ is a topological space $X$ together with a sheaf of rings $\mathcal{O}_X$ on $X$. The sheaf $\mathcal{O}_X$ is called the structure sheaf of $X$.
\end{Definition}

\begin{Definition}[Locally ringed space]
A locally ringed space is a ringed space $(X, \mathcal{O}_X)$ such that all stalks of $\mathcal{O}_X$ are local rings (i.e. they have unique maximal ideals). 
\end{Definition}

Notice: It is \underline{not} required that $\mathcal{O}_X(U)$ be a local ring for every open set $U$; in fact, this is almost never the case.        

\begin{Definition}[Affine scheme]
 Let $R$ be a commutative ring. An affine scheme is a locally ringed space $S$ such that $S \cong \Spec(R)$.
\end{Definition}

\begin{Definition}[Scheme]
 A scheme is a locally ringed space $S$ admitting a covering by open sets $U_i$, s.t. each $U_i$ (as a locally ringed space) is an affine scheme.
 
 In particular: $S$ comes with a sheaf $\mathcal{O}_S: U \mapsto \mathcal{O}_S(U)$ for every open subset $U$, the ring of regular functions on $U$.
\end{Definition}

We think of a scheme $S$ as being covered by ``coordinate charts`` which are affine schemes; we obtain a scheme $S$ by glueing together affine schemes using the Zariski topology.

\bigskip

Example: (Affine scheme)

\medskip

$\mathbb{A}_k^n=$ affine $n$-space over a field $k$, for $n \in \N$,

\smallskip

by definition we have $\mathbb{A}_k^n:= \Spec(k[x_1, \ldots, x_n])$.

\smallskip

It can be defined over any commutative ring $R$:
$\mathbb{A}_k^n:= \Spec(R[x_1, \ldots, x_n])$.

\bigskip

Example: (Scheme)

\medskip

Every affine scheme $S:= \Spec(R)$ is a scheme.

\smallskip

A polynomial $ f \in k[x_1, \ldots, x_n]$ determines a closed subscheme $ f \equiv 0$ in affine space $\mathbb{A}_k^n$, called an affine hypersurface. Formally this is $$ \Spec(k[x_1, \ldots, x_n]) / (f),$$

for example: let $k \cong \C$, then $x^2 = y^2(y+1)$ defines a singular curve in $\mathbb{A}_{\C}^2$, a nodal cubic curve.

\smallskip
$\mathbb{P}_R^n$ can be constructed as a scheme: glueing $(n+1)$ copies of affine $n$-space over $R$ along open subsets
$$ \mathbb{A}_R^n \otimes_n \mathbb{A}_R^n \otimes_n \ldots \otimes_n \mathbb{A}_R^n \; \; \; \; (n+1) \; \text{copies}$$

\begin{Remark}
An algebraic variety $X$ over a field $k$ can be defined as a scheme with certain properties. There are different conventions which schemes should be called varieties. One standard choice: a variety $X / k $ is an integral algebraic scheme of finite type over $k$.

A smooth scheme is a generalisation of the concept of a non-singular algebraic variety.

In particular: a smooth scheme of finite type over $\overline{k}$ is a non-singular algebraic variety.

Let $ k \cong \C$. A non-singular algebraic variety $X / k$ has the structure of a complex analytic manifold.
\end{Remark}

Let $S$ be a $k$-scheme.
 \begin{Definition}\label{stack}(\textbf{Algebraic stack})
  $\mathcal{X}$ is an algebraic stack if there exists a scheme $U$ and a map $U \rightarrow \mathcal{X}$ which is representable by smooth surjections.
  
  In other words: for every $S$-valued point of $\mathcal{X}$, the fiber product $U \times_{\mathcal{X}} S$ is representable by a smooth $S$-scheme $U_S$ with non-empty fibers.
 \end{Definition}

\section{Étale cohomology}

Let $S$ be a $k$-scheme and let $l$ be a prime number such that $(l,p)=1$. By $H^i(S, \Z / l^k \mathbb{Z}), \; i,k \in \N$, we define the $l$-adic cohomology group
$$H^i(S, \mathbb{Z}_l) = \underset{k}{\underset{\leftarrow}{\text{lim}}} \; H^i(S, \Z / l^k \mathbb{Z})$$ as their inverse limit.

\begin{Remark}
 The identity $\mathbb{Z}_l = \underset{k}{\underset{\leftarrow}{\text{lim}}} \; \Z / l^k \Z$ does not imply that the $l$-adic cohomology group $H^i(S, \mathbb{Z}_l)$ is equal to $H^i(S, \underset{k}{\underset{\leftarrow}{\text{lim}}} \; \Z / l^k \mathbb{Z})$ as the cohomology does not commute with taking inverse limits.
\end{Remark}

In order to remove torsion subgroups from the $l$-adic cohomology groups and to obtain cohomology groups that are vector spaces over fields of characteristic $0$, we define
$$H^i(V, \mathbb{Q}_l) = H^i(V, \mathbb{Z}_l)\otimes \mathbb{Q}_l$$

\section{Stratifications of Algebraic Stacks}

\begin{Definition}
(c-open and c-closed) A countable open (c-open) substack of an algebraic stack $\mathcal{X}$ is the complement of a countable union of Zariski closed irreducible substacks, i.e. we have $\mathcal{U}_0 = U \setminus \underset{i \in I}\cup \mathcal{A}_i$ where $\mathcal{A}_i$ are Zariski closed irreducible substacks of $U$ and $I$ is a countable index set. The complement of a c-open substack, i.e. intersections of a countable number of Zariski open substacks, is called a Zariski countable closed (c-closed) substack.
\end{Definition}

\begin{Remark}
For the existence and construction of the c-open substack $\mathcal{U}_0$ inside $U$ see \cite{PS24} for zero characteristic and \cite{SW25} for positive characteristic, where the argument is given in detail in the scheme-theoretic language.
 It uses set-theoretic ideas and in respect of Definition \ref{stack} carries over to the stack-theoretic formulation.
\end{Remark}

\begin{Definition}\cite{GL19}\label{stratification}
 Let $\mathcal{X}$ be an algebraic stack. A 
stratification of $\mathcal{X}$ consists of the following data:
 \begin{itemize}
  \item[(a)] A partially ordered set $A$.
  \item[(b)] A collection of open substacks $\{\mathcal{U}_{\alpha} \subseteq \mathcal{X}\}_{\alpha \in A}$ satisfying $\mathcal{U}_{\alpha} \subseteq \mathcal{U}_{\beta}$ when $\alpha \leq \beta$.
 \end{itemize}
This data is required to satisfy the following conditions:
\begin{itemize}
 \item[$\bullet$] For each index $\alpha \in A$, the set $\{\beta \in A: \beta \leq \alpha\}$ is finite.
 \item[$\bullet$] For every field $m$ and every map $\eta: \Spec(m) \rightarrow \mathcal{X}$, the set $$ \{\alpha \in A: \eta \text{\; factors through \;} \mathcal{U}_{\alpha}\} $$ has a least element. 
\end{itemize}
\end{Definition}

Let $\mathcal{X}$ be an algebraic stack equipped with a stratification $\{\mathcal{U}_{\alpha}\}_{\alpha \in A}$.

\begin{Notation} For each $\alpha \in A$, we let $\mathcal{X}_{\alpha}$ denote the reduced closed substack of $\mathcal{X}$ given by the complement of $\cup_{\beta < \alpha} \mathcal{U}_{\beta}$. Each $\mathcal{X}_{\alpha}$ is a locally closed substack of $\mathcal{X}$, called the strata of $\mathcal{X}$.
\end{Notation}

\begin{Remark}
 A stratification $\{\mathcal{U}_{\alpha}\}_{\alpha \in A}$ is determined by the partially ordered set $A$ together with a collection of locally closed substacks $\{\mathcal{X}_{\alpha}\}_{\alpha \in A}$: each $\mathcal{U}_{\alpha}$ is an open substack of $\mathcal{X}$ and if $m$ is a field, then a map $\eta: \Spec(m) \rightarrow \mathcal{X}$ factors through $\mathcal{U}_{\alpha}$ iff it factors through $\mathcal{X}_{\beta}$ for some $\beta \leq \alpha$. So one can identify the stratification of $\mathcal{X}$ with the locally closed substacks $\{\mathcal{X}_{\alpha} \subseteq \mathcal{X}\}_{\alpha \in A}$ (where the partial ordering of $A$ is understood to be implicitely specified).
\end{Remark}

\begin{Remark}
If $m$ is a field, then for any map $ \eta: \Spec(m) \rightarrow \mathcal{X}$ there is a 
unique index $\alpha \in A$ such that $\eta$ factors through $\mathcal{X}_{\alpha}$. In other words, $\mathcal{X}$ is the set-theoretic union of the locally closed substacks $\mathcal{X}_{\alpha}$.
\end{Remark}

\begin{Definition}{\label{convergence}}
 Let $m = \C$ or $m=\overline{\mathbb{F}_q((t))}$ and let $\mathcal{X}$ be an algebraic stack of finite type over $\Spec{(m)}$. A stratification $\{\mathcal{X}_{\alpha}\}_{\alpha \in A}$ of $\mathcal{X}$ is convergent if there exists a finite collection of algebraic stacks $\mathcal{T}_1, \ldots, \mathcal{T}_n$ over $\Spec{(m)}$  with the following properties:
 \begin{itemize}
  \item[(1)] For each $\alpha \in A$ there exists an integer $i \in \{1,2,\ldots, n\}$ and a diagram of algebraic stacks $$\mathcal{T}_i \overset{f}{\rightarrow} \overset{\sim}{\mathcal{X}}_{\alpha} \overset{g}{\rightarrow}\mathcal{X}_{\alpha}$$ where the map $f$ is a fiber bundle (locally trivial with respect to the étale topology) whose fibers are affine spaces of some fixed dimension $d_{\alpha}$ and the map $g$ is surjective, finite and radicial.
  \item[(2)] The nonnegative integers $d_{\alpha}$ in $(1)$ satisfy $\sum_{\alpha \in A}q^{-d_{\alpha}} < \infty.$
  \item[(3)] For $1 \leq i \leq n$, the algebraic stack $\mathcal{T}_i$ can be written as a stack-theoretic quotient $Y / G$, where $Y$ is an algebraic space of finite type over $m$  and $G$ is a linear algebraic group over $m$ acting on $Y$. 
 \end{itemize}
\end{Definition}

\begin{Remark}{\label{countable}}
 In definition \ref{convergence}, for $\car(m) > 0 $ the hypothesis $(2)$ guarantees that the set $A$ is at most countable.
 For $\car(m)=0$ we will assume $A$ to be at most countable. This implies that in Definition \ref{stratification}   the stratification consists of c-open substacks $\mathcal{U}_{\alpha} $.
 \end{Remark}

 We have the following
 \begin{Lemma}\label{Gysinsequence}(Gysin sequence)
Let $X$ and $Y$ be smooth quasi-projective varieties over an algebraically closed field $k$, let $g: Y \rightarrow X$ be a finite radicial morphism, and let $U \subseteq X$ be the complement of the image of $g$. Then there is a canonical fiber sequence $$C^{*-2d}(Y;\mathbb{Z}_l(-d)) \rightarrow C^*(X, \mathbb{Z}_l) \rightarrow C^*(U; \mathbb{Z}_l),$$
where $d$ denotes the relative dimension $\dim(X) - \dim(Y)$.
 \end{Lemma}        

This implies a corresponding result for algebraic stacks:
 \begin{Lemma}{\label{Gysinsequencestacks}}
  Let $\mathcal{X}$ and $\mathcal{Y}$ be smooth algebraic stacks of constant dimension over an algebraically closed field $k$, let $g: \mathcal{Y} \rightarrow \mathcal{X}$ be a finite radicial morphism, and let $\mathcal{U} \subseteq \mathcal{X}$ be the open substack of $\mathcal{X}$ complementary to the image of $g$. Then there is a canonical fiber sequence $$C^{*-2d}(\mathcal{Y};\mathbb{Z}_l(-d)) \rightarrow C^*(\mathcal{X}, \mathbb{Z}_l) \rightarrow C^*(\mathcal{U}; \mathbb{Z}_l),$$
where $d$ denotes the relative dimension $\dim(\mathcal{X}) - \dim(\mathcal{Y})$.
 \end{Lemma}
 
 \section{Construction of the stratification}
 
We now prove the main result of the paper.

\begin{Theorem}{\label{PS23.1}}
Let $U:= \Sigma \setminus \Delta$. There exists a convergent stratification $\{\mathcal{U}_{\alpha}\}_{\alpha \in A}$ of $U$ by countable open substacks for each of which the irreducibility of the monodromy representation applies in a uniform way for all $t \in \mathcal{U}_{\alpha}$ such that the monodromy argument applies in a uniform way for all $t \in U$, seen as the set-theoretic directed union $U = \underset{\underset{i}{\rightarrow}}{\cup} \; \mathcal{U}_i$.
\end{Theorem}

\begin{Proof}
 Let $d=\dim(U)$ and let $\{\mathcal{U}_{\alpha}\}_{\alpha \in A}$ be a stratification of $U$ by c-open substacks. The set $A$ is at most countable (Remark \ref{countable}). By adding additional elements to $A$ and assigning to each of those additional elements the empty substack of $U$, we may assume that $A$ is infinite.

By assumption the set $\{\beta \in A: \beta \leq \alpha\}$ is finite for each $\alpha \in A$, hence we can choose an enumeration $$ A = \{\alpha_0, \alpha_1, \alpha_2, \ldots\}$$
where each initial segment is a downward-closed subset of $A$.
 
We can then write $U$ as the union of an increasing sequence of c-open substacks $$\mathcal{U}_0 \hookrightarrow \mathcal{U}_1 \hookrightarrow \mathcal{U}_2 \hookrightarrow \ldots$$ where $\mathcal{U}_n$ is characterized by the requirement that if $m$ is a field, then a map $\eta: \Spec(m) \rightarrow \mathcal{X}$ factors through $\mathcal{U}_n$ if and only if it factors through one of the substacks $\mathcal{X}_{\alpha_0}, \mathcal{X}_{\alpha_1}, \ldots, \mathcal{X}_{\alpha_n}$.

By hypothesis, there exists a finite collection $\{\mathcal{T}_i\}_{1 \leq i\leq s}$ of smooth algebraic stacks over $\Spec(m)$, where each $\mathcal{T}_i$ has some fixed dimension $d_i$ 
and for each $n \geq 0$ there exists an index $i(n) \in \{1, \ldots, s\}$ and a diagram $$\mathcal{T}_{i(n)} \overset{f_n}{\rightarrow} \overset{\sim}{\mathcal{X}}_{\alpha_n} \overset{g_n}{\rightarrow} \mathcal{X}_{\alpha_n},$$ where $g_n$ is a finite radicial surjection and $f_n$ is an étale fiber bundle whose fibers are affine spaces of some fixed dimension $e(n)$.
The map $f_n$ induces an isomorphism in $l$-adic cohomology. Applying the Lemma \ref{Gysinsequencestacks} to the finite radicial map $$g_n: \overset{\sim}{\mathcal{X}}_{\alpha_n}\times_{\Spec(m)} \Spec(\overline{m}) \rightarrow \overline{\mathcal{U}}_n $$ we obtain fiber sequences $$C^{*-2e'_n}(\overset{\sim}{\mathcal{T}}_{i(n)};\mathbb{Z}_l(-e'_n)) \rightarrow C^*(\overline{\mathcal{U}}_n, \mathbb{Z}_l) \rightarrow C^*(\overline{\mathcal{U}}_{n-1}; \mathbb{Z}_l), $$ where $e'_n = e_n +d - d_{i(n)}$ denotes the relative dimension of the map $\overset{\sim}{\mathcal{X}}_{\alpha_n} \rightarrow U$. We have a canonical equivalence $$ \Theta:
 C^*(\overline{U};\mathbb{Z}_l) \simeq \underset{\underset{n}{\leftarrow}}{\lim}\;  C^*(\overline{\mathcal{U}}_n, \mathbb{Z}_l).$$

Passing by the fiber complexes (Gysin sequence, Lemmas \ref{Gysinsequence}, \ref{Gysinsequencestacks}) to the level of 
stacks, by the increasing sequence of c-open substacks
$$ \mathcal{U}_0 \hookrightarrow \mathcal{U}_1 \hookrightarrow \mathcal{U}_2 \hookrightarrow \ldots $$ 
we obtain $U = \underset{\underset{n}{\rightarrow}}{\cup} \; \mathcal{U}_n$.
\end{Proof}

\bigskip
\medskip

Claudia Schoemann, Simion Stoilow Institute of Mathematics of the Romanian Academy (IMAR), Bucharest, Romania, e-mail: claudia.schoemann@mathematik.uni-goettingen.de


\begin{thebibliography}{[groups]}

\bibitem[BG20]{BG20} K. Banerjee, V. Guletskii, {\em \'Etale monodromy and rational equivalence for $1$-cycles on cubic hypersurfaces
in $\mathbb{P}^5$}, Mat. Sb. Volume 211, Number 2 (2020), 3–45.

\bibitem[D74]{D74} P. \ Deligne, {\em La conjecture de Weil I}, Publ. Math\'ematiques de l'IHES {\bf} 43 (1974), 273-407.

\bibitem[DK73]{DK73} P. \ Deligne, N. Katz, {\em S\'eminaire de G\'eom\'etrie Alg\'ebrique du Bois Marie 1967-1969. Groupes de
monodromie en G\'eom\'etrie Alg\'ebrique (SGA7)} Tome II, Lecture notes in mathematics 340, Springer Verlag (1973).

\bibitem[GL19]{GL19} D. \ Gaitsgory, J. Lurie, \emph{Weil’s conjecture for Function Fields: Vol. I, II},
AMS Vol. 199 (2019).

\bibitem[La81]{La81} K. Lamotke, {\em The topology of complex projective varieties after S. Lefschetz}, Topology, Vol. 20 (1981), 15-51.

\bibitem[PS22]{PS22}
R. Paucar, C. Schoemann, {\em A theorem on zero cycles on surfaces}, Selecciones Matemáticas (2022), vol. 9, no. 1, pp 161-166, https://doi.org/10.17268/sel.mat.2022.01.13

\bibitem[PS24]{PS24} R. Paucar, C. Schoemann, {\em On the kernel of the Gysin homomorphism on Chow groups of zero cycles}, Publications Mathématiques de Bésançon (2024), pp 59-104, https://doi.org/10.5802/pmb.56

\bibitem[SW25]{SW25}
C. Schoemann, S. Werner, {\em The kernel of the Gysin homomorphism for positive characteristic}, Archiv der Mathematik, submitted (2025), https://arxiv.org/abs/2411.11417

\bibitem[Voi02]{Voi02} C. Voisin, {\em Hodge Theory and Complex Algebraic Geometry II}, Cambridge studies in advanced
mathematics 77 (2003).

\bibitem[Wei62]{Wei62} A. Weil, {\em Foundations of Algebraic Geometry}, AMS 1962.

\end{thebibliography}
\end{document}